\documentclass[11pt]{article}
\usepackage{amsmath,amssymb}
\usepackage{amsthm,bm,framed,url}
\usepackage[ruled,lined,linesnumbered]{algorithm2e}
\usepackage{graphicx}
\usepackage[monochrome]{color}
\usepackage{xcolor}
\usepackage{psfrag}
\usepackage{epstopdf}
\usepackage{wrapfig}
\usepackage{caption,subcaption}
\usepackage{cite}
\usepackage{hyperref,fixltx2e}
\usepackage{fullpage}

\newtheorem{theorem}{Theorem}

\newtheorem{lemma}{Lemma}

\DeclareMathOperator*{\minimize}{minimize}

\DeclareMathOperator*{\argmin}{arg~min}
\newcommand{\st}{\text{subject to}}
\newcommand{\half}{\frac{1}{2}}

\newcommand{\real}[1]{\Re\left\{ #1 \right\}}

\newcommand{\tr}[1]{{\rm Trace}\left\{ #1 \right\}}
\newcommand{\rank}[1]{\textsf{rank}\left( #1 \right)}
\newcommand{\inv}[1]{\left( #1 \right)^{-1}}

\newcommand{\R}{\mathcal{R}}
\newcommand{\C}{\mathcal{C}}
\renewcommand{\H}{\mathcal{H}}
\renewcommand{\O}{\mathcal{O}}

\newcommand{\N}{\mathcal{N}}
\newcommand{\CN}{\mathcal{CN}}

\newcommand{\x}{\bm{x}}
\newcommand{\y}{\bm{y}}
\newcommand{\z}{\bm{z}}
\newcommand{\zt}{\tilde{\z}}
\renewcommand{\u}{\bm{u}}

\newcommand{\s}{\bm{s}}

\newcommand{\eye}{\bm{I}}

\newcommand{\X}{\bm{X}}

\newcommand{\A}{\bm{A}}
\renewcommand{\a}{\bm{a}}
\renewcommand{\b}{\bm{b}}
\newcommand{\Q}{\bm{Q}}
\renewcommand{\L}{\bm{\varLambda}}
\newcommand{\bt}{\tilde{\bm{b}}}
\newcommand{\tb}{\tilde{b}}

\newcommand{\bzeta}{\bm{\zeta}}
\newcommand{\tzeta}{\tilde{\zeta}}
\newcommand{\btzeta}{\tilde{\bm{\zeta}}}

\newcommand{\w}{\bm{w}}
\newcommand{\h}{\bm{h}}
\newcommand{\g}{\bm{g}}

\providecommand{\keywords}[1]{\textbf{Keywords:} #1}

\newcommand{\reminder}[1]{} 

\begin{document}

\title{\bf Consensus-ADMM for General Quadratically Constrained Quadratic Programming}
\author{Kejun~Huang
and
Nicholas~D.~Sidiropoulos
}
\maketitle

\begin{abstract}
Non-convex quadratically constrained quadratic programming (QCQP) problems have numerous applications in signal processing, machine learning, and wireless communications, albeit the general QCQP is NP-hard, and several interesting special cases are NP-hard as well. This paper proposes a new algorithm for general QCQP. The problem is first reformulated in consensus optimization form, to which the alternating direction method of multipliers (ADMM) can be applied. The reformulation is done in such a way that each of the sub-problems is a QCQP with only one constraint (QCQP-1), which is efficiently solvable irrespective of (non-)convexity. The core components are carefully designed to make the overall algorithm more scalable, including efficient methods for solving QCQP-1, memory efficient implementation, parallel/distributed implementation, and smart initialization. The proposed algorithm is then tested in two applications: multicast beamforming and phase retrieval. The results indicate superior performance over prior state-of-the-art methods.
\end{abstract}

\section{Introduction}
Quadratically constrained quadratic programming (QCQP) is an optimization problem that minimizes a quadratic function subject to quadractic inequality and equality constraints~\cite{luo2010semidefinite}. We write it in the most general form as follows:
\begin{equation}\label{prob:qcqp}
\begin{aligned}
\minimize_{\x\in\C^n}~~ & \x^H\A_0\x - 2\real{\b_0^H\x}, \\
\st~~	& \x^H\A_i\x - 2\real{\b_i^H\x} \leq c_i,\\
		& ~~~~~~~~~\forall~i=1,...,m.
\end{aligned}
\end{equation}
Notice that for simplicity we only write the constraints as inequalities, but they can be equalities as well (each can be expressed as two inequalities).

A QCQP is in general NP-hard, except for some special cases, for example when all the $\{\A_i\}_{i=1}^m$ in ``$\leq$'' inequality constraints are positive semi-definite~\cite[\S~4.4]{boyd2004convex}, $m$ is ``small''~\cite{ye2003new,beck2006strong,huang2010rank,ai2011new}, or if the quadratic terms are all homogeneous and $\{\A_i\}_{i=0}^m$ are all Toeplitz~\cite{konar2015hidden}.

For general non-convex QCQPs, the prevailing method to tackle the problem is through semi-definite relaxation (SDR), where the following semi-definite programming (SDP) problem is solved instead
\begin{equation}\label{prob:sdr}
\begin{aligned}
\minimize_{\X\in\H^n,\x\in\C^n}~~ & \tr{\A_0\X} - 2\real{\b_0^H\x}, \\
\st~~	& \begin{bmatrix}
		\X & \x \\
		\x^H & 1
		\end{bmatrix} \succeq 0, \\
		& \tr{\A_i\X} - 2\real{\b_i^H\x} \leq c_i,\\
		& ~~~~~~~~~\forall~i=1,...,m,
\end{aligned}
\end{equation}
where $\H^n$ denotes the set of $n \times n$ complex Hermitian matrices. Problem (\ref{prob:sdr}) is obtained by relaxing the otherwise equivalent constraint $\X=\x\x^H$ to a convex one $\X \succeq \x\x^H$. After solving the convex relaxation problem (\ref{prob:sdr}), we not only get a non-trivial lower bound on the optimal cost of (\ref{prob:qcqp}), but also possibly a solution of (\ref{prob:qcqp}) if the solution $\X_\star$ of (\ref{prob:sdr}) turns out to be rank one. If this is not the case, in certain cases there is an efficient way to generate approximate solutions to the original problem in (\ref{prob:qcqp}) from the higher-rank solution of (\ref{prob:sdr}). Let $(\x_\star,\X_\star)$ be an optimal solution of (\ref{prob:sdr}), then by drawing random points $\x\sim\CN(\x_\star,\X_\star)$, possibly followed by a simple projection or scaling if applicable, one can obtain approximate solutions not far away from the SDR lower bound, for certain kinds of QCQP problems. That is, in certain cases it has been shown that this randomization step (with sufficient draws) is guaranteed to yield some quantified sub-optimality, see for example~\cite{goemans1995improved,luo2007approximation}.

If the problem dimension in (\ref{prob:qcqp}) is large, then squaring the number of variables as in (\ref{prob:sdr}) makes the latter very difficult to solve. If a general purpose SDP solver is used to solve (\ref{prob:sdr}) using the interior point method, the worst case complexity can be as high as $\O(n^{6.5})$. Another downside of SDR is that, if it is not obvious how to obtain a feasible point for the constraint set of (\ref{prob:qcqp}), in a lot of cases SDR randomization will not give us a feasible point either.

Another way to tackle problem (\ref{prob:qcqp}) is through convex restriction, also known as successive convex approximation (SCA) or convex-concave procedure (CCP) for the more general difference of convex programming (DCP) problem~\cite{yuille2003concave}. Noticing that any Hermitian matrix can be written as the difference of two positive semi-definite matrices, we can denote
\[
\A_i = \A_i^{(+)} + \A_i^{(-)},~\forall~i=0,1,...,m,
\]
where $\A_i^{(+)} \succeq 0$ and $\A_i^{(-)} \preceq 0$. Then for each quadratic term, we have that
\[
\x^H\A_i\x \leq \x^H\A_i^{(+)}\x + 2\real{\x^{(0)H}\A_i^{(-)}\x} - \x^{(0)H}\A_i^{(-)}\x^{(0)},
\]
for any point $\x^{(0)}\in\C^n$. Therefore, starting with an initial point $\x^{(0)}$, we can iteratively solve the following convex QCQP until we obtain an approximate solution of (\ref{prob:qcqp})
\begin{align*}
\x^{(t+1)} \leftarrow & \\
\arg\min_{\x}~~ & \x^H\A_0^{(+)}\x + 2\real{\x^{(t)H}\A_0^{(-)}\x} - 2\real{\b_0^H\x} \\
\st~~ & \x^H\A_i^{(+)}\x + 2\real{\x^{(t)H}\A_i^{(-)}\x} - 2\real{\b_i^H\x} \\
	& \leq c_i + \x^{(t)H}\A_i^{(-)}\x^{(t)}, \forall~i=1,...,m.
\end{align*}
This is a convex restriction because each quadratic function is replaced with its convex upper bound function. If we start with a feasible $\x^{(0)}$, then it is easy to show that the sequence $\{\x^{(t)}\}$ will remain feasible afterwards. However, if we start with an infeasible $\x^{(0)}$, it is possible (and often the case) that the restricted constraint set becomes empty, thus the iterates cannot proceed. Recently \cite{mehanna2015feasible} proposed {\em feasible point pursuit - successive convex approximation} (FPP-SCA) to address this issue, by adding a slack variable to each constraint and an $\ell_1$ penalty on the slacks to the cost. FPP-SCA produces good results in both finding a feasible point and approaching closer to the SDR lower bound. The potential disadvantage of FPP-SCA is that in each iteration we still need to solve a non-trivial convex optimization problem, which may take a lot of time even for a moderate number of iterations, if a general-purpose convex optimization solver is used for large-scale problems.

In this paper, we propose a rather different approach to handle QCQPs. The proposed algorithm is based on two building blocks:
\begin{enumerate}
\item Any QCQP with only one constraint (QCQP-1) can be solved to optimality, and in various cases this can be done efficiently;
\item Adopting the alternating direction method of multipliers (ADMM) for consensus optimization as the general algorithmic framework, problem (\ref{prob:qcqp}) can decomposed into $m$ QCQP-1's in each iteration, thus leading to efficient optimal updates.
\end{enumerate}

In the rest of this section, we briefly review the basics of the two aforementioned building blocks. Then the general algorithmic framework is introduced in Section~\ref{sec:2}. In Section~\ref{sec:3}, we look into one specific step of the algorithm, and explain how this seemingly non-trivial QCQP-1 sub-problem can be solved very efficiently. Some detailed implementation issues are described in Section~\ref{sec:4}, including a memory-efficient implementation for certain types of constraints, an empirical way of initialization that works very well in practice, and discussions on parallel and distributed implementations with small communication overhead. Simulation results are presented in Section~\ref{sec:5}, where the proposed algorithm is used for feasible point pursuit, multicast beamforming, and phase retrieval, showing great performance and versatility in various kinds of non-convex QCQP problems. Conclusions are drawn in Section~\ref{sec:6}.

\subsection{QCQP with only one constraint (QCQP-1)}
One of the most fundamental results in QCQP is that, any QCQP with only one constraint (QCQP-1) can be solved optimally, despite the fact that the quadratic terms may be indefinite. The fundamental idea behind this result is the following lemma~\cite[Appendix~B]{boyd2004convex}.
\begin{lemma}
For all $\X,\A,\bm{B}\in\H^n$, and $\X\succeq 0$, there exists an $\x\in\C^n$ such that
\[
\x^H\A\x = \tr{\A\X},~\x^H\bm{B}\x = \tr{\bm{B}\X}.
\]
\end{lemma}
This means that after we find a solution for the SDR of a QCQP-1, regardless of its rank, we can always find an equivalent rank one solution.
{\color{blue}Our experience from simulations is that SDR seems to always return a rank one solution for a QCQP-1. Even if not, one can resort to rank reduction as in \cite{huang2010rank}, which handles a more general rank reduction problem.}
This result is also closely related to the {\em generalized eigenvalue problem} in linear algebra and the {\em S-procedure} in control.

\subsection{Consensus optimization using ADMM}
Now we briefly introduce the algorithmic tool to be used in this paper, which is based on the alternating direction method of multipliers (ADMM) \cite{boyd2011distributed}. Consider the following optimization problem
\[
\minimize_{\x}~~\sum_{i=1}^{m}f_i(\x) + r(\x),
\]
in which the main objective is to minimize a sum of cost functions $f_1$, ..., $f_m$, subject to some additional regularization $r$ on $\x$. To solve it using ADMM, we first reformulate it into a \emph{consensus} form by introducing $m$ auxiliary variables $\z_1$, ..., $\z_m$, as
\begin{align*}
\minimize_{\x,\{\z_i\}_{i=1}^m}~~ & \sum_{i=1}^{m}f_i(\z_i) + r(\x), \\
\st~~ & \z_i=\x,~\forall i=1,...,m.
\end{align*}
Then we can easily write down the (scaled-form) ADMM iterates for this problem as
\begin{align*}
\x\; & \leftarrow \arg\min_{\x}~ r(\x) + \rho\sum_{i=1}^{m}\|\z_i-\x+\u_i\|^2, \\
\z_i & \leftarrow \arg\min_{\z_i}~f_i(\z_i) + \rho\|\z_i-\x+\u_i\|^2,~\forall~i=1,...,m,\\
\u_i & \leftarrow \u_i + \z_i - \x,~\forall~i=1,...,m,
\end{align*}
where $\x$ is treated as the first block, the set of auxiliary variables $\{\z_i\}$ are treated as the second block, and $\u_i$ is the scaled dual variable corresponding to the equality constraint $\z_i=\x$.

There are several advantages of this consensus-ADMM algorithm. First and foremost, it is designed for distributed optimization, since each $\z_i$ can be updated in parallel; through careful splitting of the sum of the cost functions, we can also make each update very efficient (possibly in closed-form), which may not be the case for the batch problem. Finally, since it falls into the general algorithmic framework of ADMM, it converges as long as the problem is convex, for all $\rho>0$.

\section{General algorithmic framework}\label{sec:2}
We now describe how to apply consensus-ADMM for general QCQPs. Let us first transform (\ref{prob:qcqp}) into a consensus form
\begin{equation}\label{prob:qcqp-c}
\begin{aligned}
\minimize_{\x,\{\z_i\}_{i=1}^m}~~ & \x^H\A_0\x - 2\real{\b_0^H\x}, \\
\st~~	& \z_i^H\A_i\z_i - 2\real{\b_i^H\z_i} \leq c_i, \\
		& \z_i = \x,~\forall~i=1,...,m,
\end{aligned}
\end{equation}
then the corresponding consensus-ADMM algorithm takes the form of the following iterations:
\begin{equation}\label{alg:admm}
\begin{aligned}
\x\; & \leftarrow \inv{\A_0+m\rho\eye}\left(\b_0 + \rho\sum_{i=1}^{m} \left(\z_i + \u_i\right)\right), \\
\z_i\, & \leftarrow \arg\min_{\z_i}~\|\z_i-\x+\u_i\|^2, \\
		& ~~~\st~~ \z_i^H\A_i\z_i - 2\real{\b_i^H\z_i} \leq c_i \\
\u_i & \leftarrow \u_i + \z_i - \x.
\end{aligned}
\end{equation}

The reason we put our algorithm into this form is based on the fact that each update of $\z_i$ is a QCQP-1, thus we know it can be updated optimally, despite the fact that the quadratics may be indefinite. The update for $\x$ is an unconstrained quadratic minimization, and for an indefinite $\A_0$ we need to choose a large enough $\rho$ to ensure that the minimum is not unbounded; if $\A_0 + \rho m\eye \succeq 0$ is satisfied then the solution is simply given by solving a linear equation, and we can cache the Cholesky factorization of $\A_0 + m \rho \eye$ to save computations in the subsequent iterations.

\subsection{Convergence}
ADMM was first designed for convex problems, for which it is known to converge under mild conditions~\cite{eckstein1992douglas}. Despite the lack of theoretical guarantees, ADMM has also been used for non-convex problems, see for example~\cite[\S9]{boyd2011distributed}. In~\cite{xu2012alternating}, ADMM was applied to non-negative matrix factorization (a non-convex problem) with missing values, and it was shown that, {\em if} ADMM converges for this non-convex problem, {\em then} it converges to a KKT point. Some follow-up works on other non-convex problems have made similar claims~\cite{jiang2014alternating,liavas2015parallel}. A stronger result was recently shown in~\cite{hong2014convergence}, where it was proven that (there exists a convergent subsequence and) {\em every limit point is a stationary point} for a class of non-convex consensus and sharing problems. The proof in \cite{hong2014convergence} assumes Lipschitz continuity of the non-convex cost functions (not constraints) to establish that the augmented Lagrangian function is non-increasing, provided the parameter $\rho$ is large enough.

Unfortunately, the convergence result in~\cite{hong2014convergence} cannot be applied to our algorithm here, even though both are  dealing with non-convex \emph{consensus} problems. The very first step in the proof of \cite{hong2014convergence} shows that the augmented Lagrangian is monotonically non-increasing under certain conditions. These conditions include Lipschitz continuity of the non-convex cost functions (but not the constraints) and that the parameter $\rho$ is large enough. If we want to borrow the arguments made in \cite{hong2014convergence}, we would need to first establish the monotonicity of the augmented Lagrangian. However, our numerical experience is that the augmented Lagrangian is not monotonic, even if we set $\rho$ to be very large.
Therefore, we limit ourselves to the following weaker convergence result.

{\color{blue}
\begin{theorem}\label{thm:convergence}
Denote $\x^t$ and $\z_i^t$ the updates obtained at the $t$-th iteration of Algorithm~\eqref{alg:admm}. Assume that the $\z_i^t$'s are well-defined for all $t$ and $i$, and that
\[
\lim_{t\rightarrow+\infty}(\z_i^t - \x^t) = 0, \forall i=1,...,m,
\]
and
\[
\lim_{t\rightarrow+\infty}(\x^{t+1} - \x^t) = 0,
\]
then any limit point of $\{\x^t\}$ is a KKT point of \eqref{prob:qcqp}
\end{theorem}
}
\begin{proof}
See Appendix~\ref{Appendix:proof}.
\end{proof}

\section{Efficient \texorpdfstring{$\z_i$}{z\_i}-updates}\label{sec:3}
Now let us focus on the update of $\z_i$. From our previous discussion on QCQP-1 we know that the update of $\z_i$ can always be solved to optimality, by strong duality; in other words, if we solve the SDR of a QCQP-1, we are guaranteed to obtain a rank one solution. However, with a number of $\z_i$ to be updated iteratively, it is not desirable to rely on general SDP algorithms to update $\z_i$. Therefore, we now take a detailed look into QCQP-1, and show how to solve it efficiently. For ease of notation, let us drop the subscript, define $\bzeta=\x-\u$, and denote the sub-problem as
\begin{equation}\label{prob:qcqp1}
\begin{aligned}
\minimize_{\z}~~& \left\|\z-\bzeta\right\|^2 \\
\st~~ & \z^H\A\z - 2\real{\b^H\z} = c
\end{aligned}
\end{equation}
We changed the constraint to equality here to simplify subsequent derivations. For an inequality constraint, we first check whether $\bzeta$ is feasible: if yes, then $\bzeta$ is the solution; if not, then the constraint must be satisfied as equality, according to complementary slackness, thus the following method can be applied.

We start from simpler cases, and gradually build up to the most general case.

\subsection{\texorpdfstring{$\rank{\A}=1$}{rank(A)=1}, \texorpdfstring{$\b = {\bf 0}$}{b=0}}\label{sec:qcqp1-1}
For this simple case, the constraint can be equivalently written as
\[
|\a^H\z|^2=c,
\]
or simply as a linear constraint with an unknown phase
\[
\a^H\z = \sqrt{c} e^{j\theta}.
\]
Assuming we know $\theta$, problem (\ref{prob:qcqp1}) becomes a simple projection onto an affine subspace, for which we know the solution is given by
\[
\z = \bzeta + \left( \frac{\sqrt{c} e^{j\theta}-\a^H\bzeta}{\|\a\|^2} \right)\a.
\]
Plugging this back to the objective, it is easy to see that the minimum is attained if we choose $\theta$ to be the angle of $\a^H\bzeta$. Therefore, the update of $\z$ in this case is given by
\begin{equation}\label{eq:A1b0}
\z = \bzeta+ \frac{\sqrt{c}-|\a^H\bzeta|}{\|\a\|^2|\a^H\bzeta|}\a\a^H\bzeta.
\end{equation}

In the real case, the unknown phase becomes an unknown sign, and similar steps can be made to result in the same closed-form solution (\ref{eq:A1b0}).

\subsection{\texorpdfstring{$\rank{\A}>1$}{rank(A)>1}, \texorpdfstring{$\b = {\bf 0}$}{b=0}}\label{sec:qcqp1-3}
For $\A$ with a higher rank, there is in general no closed-form solution for (\ref{prob:qcqp1}). However, it is still possible to efficiently update $\z$. Let the eigen-decomposition of $\A$ be $\Q\L\Q^H$, where $\L$ is diagonal real and $\Q$ is unitary, because $\A$ is Hermitian. Define $\zt=\Q^H\z$, $\btzeta=\Q^H\bzeta$, then the problem is equivalent to
\begin{align*}
\minimize_{\zt}~~ & \|\zt-\btzeta\|^2, \\
\st~~ & \zt^H\L\zt = c.
\end{align*}
The corresponding Lagrangian is
\[
L = \|\zt-\btzeta\|^2 + \mu\left(\zt^H\L\zt - c\right),
\]
with a single Lagrange multiplier $\mu$. A necessary condition for optimality is that $\nabla L = 0$, i.e.,
\begin{align*}
\nabla L = 2(\zt-\btzeta) + 2\mu\L\zt = 0,
\end{align*}
therefore, $\zt = \inv{\eye+\mu\L}\btzeta$. Plugging this solution back into the equality constraint, we have
\[
\btzeta^H\inv{\eye+\mu\L}\L\inv{\eye+\mu\L}\btzeta=c,
\]
or equivalently
\[
\sum_{k=1}^{n} \frac{\lambda_k}{(1+\mu\lambda_k)^2}|\tzeta_k|^2 = c,
\]
which means the correct Lagrange multiplier $\mu$ can be numerically found by solving this nonlinear equation, via for example bisection or Newton's method. In fact, we can also show that the desired solution is unique, leaving no ambiguity to the value of $\mu$. From the dual of a \mbox{QCQP-1}, we have that $\eye+\mu\L\succeq0$~\cite[Appendix~B]{boyd2004convex}, i.e.,
\[
1+\mu\lambda_k \geq 0,~\forall~k=1,...,n.
\]
This can give us a first possible region where the correct $\mu$ can be: $\mu\leq-1/\lambda_{\rm min}$ if $\lambda_{\rm min}<0$, and $\mu\geq-1/\lambda_{\rm max}$ if $\lambda_{\rm max}>0$.

Moreover, if we define
\[
\phi(\mu) = \sum_{k=1}^{n} \frac{\lambda_k}{(1+\mu\lambda_k)^2}|\tzeta_k|^2 - c,
\]
then
\[
\phi'(\mu) = -2\sum_{k=1}^{n} \frac{\lambda_k^2}{(1+\mu\lambda_k)^3}|\tzeta_k|^2,
\]
and for all $\mu$ such that $\eye+\mu\L\succeq0$, $\phi'(\mu) < 0$, which means $\phi(\mu)$ is monotonically decreasing (strictly) in that region, therefore the solution for $\phi(\mu)=0$ is unique.
{\color{blue}In fact, we can show that there exists a root within that interval, as long as the constraint set is not empty: if $-1/\lambda_{\rm max}\leq \mu \leq -1/\lambda_{\rm min}$, then $\phi(-1/\lambda_{\rm max})=+\infty$, and $\phi(-1/\lambda_{\rm min})=-\infty$, which together with the monotonicity imply that a root always exists in between. If $\bm{\varLambda}\succeq0$, the interval becomes $-1/\lambda_{\rm max}\leq \mu \leq +\infty$, consequently $-c\leq\phi(\mu)\leq+\infty$, so a root exists if and only if $c\geq0$, but if $c<0$ then $\tilde{\z}^T\bm{\varLambda}\tilde{\z}=c$ is infeasible. A similar argument applies to the case when $\bm{\varLambda}\preceq0$.}
Once the value of $\mu$ is found, we can plug it back to obtain $\zt$, and the desired update of $\z$ is simply given by $\z=\Q\zt$.

To save computation, we can cache the eigen-decomposition of $\A$. Then in the subsequent ADMM iterations the computation is dominated by the matrix-vector multiplication $\Q\zt$, since evaluating either $\phi(\mu)$ or $\phi'(\mu)$ (if Newton's method is used) only takes $O(n)$ complexity.

\subsection{\texorpdfstring{$\rank{\A}>1$}{rank(A)>1}, \texorpdfstring{$\b \neq {\bf 0}$}{b not zero}}\label{sec:qcqp1-4}
Now we have reached the most general case when $\A$ can have higher rank and $\b$ can be non-zero. The idea is very similar to the previous case, although the expressions are a little more complicated. Again let $\A=\Q\L\Q^H$ be the eigen-decomposition, problem (\ref{prob:qcqp1}) is equivalent to
\begin{align*}
\minimize_{\zt}~~ & \|\zt-\btzeta\|^2, \\
\st~~ & \zt^H\L\zt - 2\real{\bt^H\zt} = c,
\end{align*}
Where $\zt=\Q^H\z$, $\btzeta=\Q^H\bzeta$, and $\bt = \Q^H\b$. Setting the gradient of the Lagrangian equal to zero, we have
\[
\zt = \inv{\eye+\mu\L}(\btzeta+\mu\bt).
\]
Plugging it back to the equality constraint, it becomes a nonlinear equation with respect to $\mu$,
\[
\phi(\mu)=\sum_{k=1}^{n}\lambda_k\left|\frac{\tzeta_k+\mu\tb_k}{1+\mu\lambda_k}\right|^2 -
			2\real{ \sum_{k=1}^{n}\tb_k^*\frac{\tzeta_k+\mu\tb_k}{1+\mu\lambda_k} } - c,
\]
and its derivative
\[
\phi'(\mu) = -2\sum_{k=1}^{n}\frac{|\tb_k-\lambda_k\tzeta_k|^2}{(1+\mu\lambda_k)^3} < 0,
\]
for all $\mu$ such that $\eye+\mu\L\succeq 0$, which is necessary for optimality of (\ref{prob:qcqp1}). Therefore, $\phi(\mu)$ is monotonic in the possible region of solution, and any local solution (for example found by bisection or Newton's method) is guaranteed to be the unique (thus correct) solution, {\color{blue}which always exists for a non-empty constraint set, similar to the previous case}. Notice that if $\b={\bf 0}$, $\phi(\mu)$ and $\phi'(\mu)$ reduce to the simpler expression that we derived in the previous subsection. Detailed implementation of bisection and Newton's method to solve $\phi(\mu)=0$ is given in Alg.~\ref{algo:bisec}~and~\ref{algo:newton}. In practice, bisection converges linearly ($\sim$ 20 iterations) while Newton's method converges quadratically ($\sim$ 5 iterations), but bisection is numerically more stable, so the best choice is application-specific.

\begin{algorithm}[t]
$\lceil\mu\rceil   \leftarrow \text{max number available in float representation}$\;
$\lfloor\mu\rfloor \leftarrow \text{min number available in float representation}$\;
\lIf{\upshape $\lambda_\text{min}<0$}{$\lceil\mu\rceil \leftarrow -1/\lambda_\text{min}$}
\lIf{\upshape $\lambda_\text{max}>0$}{$\lfloor\mu\rfloor \leftarrow -1/\lambda_\text{max}$}
\Repeat{$\lceil\mu\rceil-\lfloor\mu\rfloor<\varepsilon$}{
$\mu = (\lceil\mu\rceil+\lfloor\mu\rfloor)/2$\;
\lIf{$\phi(\mu)>0$}{$\lfloor\mu\rfloor\leftarrow\mu$}
\lElse{$\lceil\mu\rceil\leftarrow\mu$}
}
\Return{$\mu = (\lceil\mu\rceil+\lfloor\mu\rfloor)/2$}
\caption{solving $\phi(\mu)=0$ using bisection}
\label{algo:bisec}
\end{algorithm}
\begin{algorithm}[t]
$\mu \leftarrow {\color{blue}-}(\lambda_\text{min}+\lambda_\text{max})/2\lambda_\text{min}\lambda_\text{max}$\;
\Repeat{$-\phi(\mu)^2/\phi'(\mu)<\varepsilon$}{
$\mu = \mu - \phi(\mu)/\phi'(\mu)$\;
}
\Return{$\mu$}
\caption{solving $\phi(\mu)=0$ using Newton's method}
\label{algo:newton}
\end{algorithm}

An interesting observation from this most general case is that, solving a QCQP-1 always boils down to solving a scalar nonlinear equation $\phi(\mu)=0$. It is easy to see that if $\A$ has $p$ distinct eigenvalues, solving $\phi(\mu)=0$ is equivalent to solving a polynomial of degree $2p+1$ ($2p$ if $\b={\bf 0}$). Polynomials of order $\geq 5$ do not admit closed-form expressions for their roots, necessitating the use of numerical methods like bisection or Newton's method. 

{\color{blue}
\paragraph*{Remark}
So far we have assumed that $\eye+\mu\L$ is invertible, which may not always be the case. However, recall that the duality of QCQP-1 implies $\eye+\mu\L\succeq0$, therefore there are at most two possible values of $\mu$ that can make the matrix $\eye+\mu\L$ singular:
$\mu=-1/\lambda_{\rm min}$ if $\lambda_{\rm min}<0$, and $\mu=-1/\lambda_{\rm max}$ if $\lambda_{\rm max}>0$, so for completeness one may first check these two values of $\mu$, although this situation never occurred in our experiments. 
}

\subsection{Bound constraint}\label{sec:qcqp1-5}
The basic idea of making the $\z$-updates equivalent to solving a QCQP-1 is that the latter is always efficiently solvable. In some cases this efficiency can be maintained even if we incorporate some more constraints. One such case is that of a quadratic term that is bounded from both sides, i.e.,
\begin{align*}
\minimize_{\z}~~& \left\|\z-\bzeta\right\|^2, \\
\st~~ & c-\epsilon \leq \z^H\A\z - 2\real{\b^H\z} \leq c+\epsilon.
\end{align*}
Using the same idea as before, we can write down the Lagrangian and the solution again takes the form
\[
\z = \Q\inv{\eye+\mu\L}\Q^H(\bzeta+\mu\b),
\]
where $\A=\Q\L\Q^H$ is the eigen-decomposition of $\A$, and $\mu$ is such that
\[
\begin{cases}
\mu = 0, & \quad \text{if } -\epsilon\leq\phi(0)\leq\epsilon, \\
\mu < 0, & \quad \text{then } \phi(\mu) =-\epsilon, \\
\mu > 0, & \quad \text{then } \phi(\mu) = \epsilon.
\end{cases}
\]
In fact, since we know $\phi(\mu)$ is monotonically decreasing within the feasible region, if $\phi(0)>\epsilon$, both the solution of $\phi(\mu)=\pm\epsilon$ are positive, therefore we must take the solution of $\phi(\mu)=\epsilon$, and vice versa. In other words, if
\[
c-\epsilon \leq \bzeta^H\A\bzeta - 2\real{\b^H\bzeta} \leq c+\epsilon,
\]
then $\z=\bzeta$; if it is greater than $c+\epsilon$, the upper-bound constraint must be active, and like-wise if it is less than $c-\epsilon$. This is very intuitive, since we are just ``rounding'' the constraint to the closest bound.

\section{Implementation issues}\label{sec:4}
So far we have derived an ADMM algorithm for general QCQP problems, which features straightforward iterations and efficient per-iteration updates. In this section we revisit the entire algorithm and discuss detailed implementations to make it more actionable.

\subsection{Memory-efficient implementation}
An apparent disadvantage of our algorithm is that we need to introduce an auxiliary variable $\z_i$ and the corresponding dual variable $\u_i$ for every single quadratic constraint. For $\x\in\C^n$ and $m$ such constraints, we need $\O(mn)$ memory just to store the intermediate variables. Depending on the application, this memory requirement may be too demanding. For example, if $\A_i=\a_i^{}\a_i^{H},~\forall~i=1,...,m$, it only takes $\O(mn)$ memory to describe the problem, or even as small as $\O(m)$ if the $\a_i$'s are highly structured, e.g., obtained from the rows of the discrete Fourier transform (DFT) matrix. In such cases $\O(mn)$ intermediate memory seems very unappealing for large $m$ and $n$. This is less of an issue when the $\A_i$'s are all full rank, since then we need $\O(mn^2)$ memory to just specify the problem, and if that is affordable, then $\O(mn)$ memory for intermediate variables seems relatively reasonable.

Consider the following special QCQP, which occurs frequently in practice:
\begin{equation}\label{prob:rank_1}
\begin{aligned}
\minimize_{\x}~~ & \x^H\A_0\x - 2\real{\b_0^H\x} \\
\st~~ & |\a_i^H\x|^2=c_i,~\forall~i=1,...,m.
\end{aligned}
\end{equation}
Again, the algorithm that follows can be easily modified to tackle inequality constraints or bound constraints, but we start with equality constraints here for clarity.
According to our previous discussion, we can write down explicitly the consensus-ADMM iterations as
\begin{equation}\label{algo:rank_1}
\begin{aligned}
\x   & \leftarrow \inv{\A_0+m\rho\eye}\left(\b_0 + \rho\sum_{i=1}^{m} \left(\z_i + \u_i\right)\right), \\
\z_i & \leftarrow \x-\u_i + \frac{\sqrt{c_i}-|\a_i^H(\x-\u_i)|}{\|\a_i\|^2|\a_i^H(\x-\u_i)|}\a_i\a_i^H(\x-\u_i), \\
\u_i & \leftarrow \u_i + \z_i - \x.
\end{aligned}
\end{equation}
Define $\z_s = \sum_i \z_i$ and $\u_s = \sum_i \u_i$, then we can simplify the algorithm as
\begin{align*}
\x   & \leftarrow \inv{\A_0+m\rho\eye}\left(\b_0 + \z_s + \u_s\right), \\
\z_s & \leftarrow m\x - \u_s + \A_s\bm{\nu}, \\
\u_s & \leftarrow \u_s + \z_s - m\x,
\end{align*}
where $\A_s=[~\a_1~\a_2~...~\a_m~]$ is a $n \times m$ matrix formed by parallel stacking all the $\a_i$ vectors as its columns, and $\bm{\nu}$ is a vector of length $m$ with its elements defined as
\[
\nu_i = \frac{\sqrt{c_i}-|\a_i^H(\x-\u_i)|}{\|\a_i\|^2|\a_i^H(\x-\u_i)|}\a_i^H(\x-\u_i).
\]
If we are given the vector $\bm{\nu}$ at every iteration, then we can simply work with the summation of the local variables $\z_s$ and $\u_s$ without the possible memory explosion. To compute the vector $\bm{\nu}$, we notice that it is sufficient to know the value of $\a_i^H(\x-\u_i)$ for each $\nu_i$, and since we keep track of $\x$ explicitly, the only difficulty is to keep track of $\a_i^H\u_i$ without the actual value of $\u_i$. At the end of each iteration, by combining the updates of $\z_i$ and $\u_i$, we have that
\[
\u_i \leftarrow \frac{\sqrt{c_i}-|\a_i^H(\x-\u_i)|}{\|\a_i\|^2|\a_i^H(\x-\u_i)|}\a_i^{}\a_i^H(\x-\u_i),
\]
therefore
\[
\a_i^H\u_i \leftarrow \frac{\a_i^H(\x-\u_i)}{|\a_i^H(\x-\u_i)|}\left(\sqrt{c_i}-|\a_i^H(\x-\u_i)|\right).
\]
Now if we define $\alpha_i=\a_i^H\u_i$, it is apparent that we can update $\alpha_i$ iteratively as
\[
\alpha_i \leftarrow \frac{\a_i^H\x-\alpha_i}{|\a_i^H\x-\alpha_i|}\left(\sqrt{c_i}-|\a_i^H\x-\alpha_i|\right).
\]

To sum up, a memory-efficient way to implement consensus-ADMM for problem (\ref{prob:rank_1}) takes the form
\begin{equation}\label{algo:efficient}
\begin{aligned}
\x   & \leftarrow \inv{\A_0+m\rho\eye}\left(\b_0 + \rho(\z_s + \u_s)\right), \\
\bm{\xi} & \leftarrow \A_s^H\x, \\
\nu_i  & \leftarrow \frac{\xi_i-\alpha_i}{|\xi_i-\alpha_i|}\frac{\sqrt{c_i}-|\xi_i-\alpha_i|}{\|\a_i\|^2}, \\
\z_s & \leftarrow m\x - \u_s + \A_s\bm{\nu}, \\
\u_s & \leftarrow \u_s + \z_s - m\x, \\
\alpha_i & \leftarrow \frac{\xi_i-\alpha_i}{|\xi_i-\alpha_i|}\left(\sqrt{c_i}-|\xi_i-\alpha_i|\right).
\end{aligned}
\end{equation}
The explicit variables are $\x,\z_s,\u_s\in\C^n$ and $\bm{\xi},\bm{\nu},\bm{\alpha}\in\C^m$, so the total memory consumption now is $\O(m+n)$, compared to $\O(mn)$ in the original form.

Finally, we show that the modified iterates in (\ref{algo:efficient}) can handle some variations in the constraints. Suppose the $i$-th constraint is an inequality $|\a_i^H\x|^2\leq c_i$, then the update of $\z_i$ in (\ref{algo:rank_1}) should be
\[
\z_i \leftarrow \left\{
\begin{array}{l}
\x-\u_i , \text{~if~} |\a_i^H(\x-\u_i)|^2\leq c_i \\
\x-\u_i + \frac{\sqrt{c_i}-|\a_i^H(\x-\u_i)|}{\|\a_i\|^2|\a_i^H(\x-\u_i)|}\a_i\a_i^H(\x-\u_i), \text{ otherwise,}
\end{array}
\right.
\]
or simply
\[
\z_i \leftarrow \x-\u_i +
	\frac{\Big[\sqrt{c_i}-|\a_i^H(\x-\u_i)|\Big]_-}{\|\a_i\|^2|\a_i^H(\x-\u_i)|}\a_i\a_i^H(\x-\u_i).
\]
This means the corresponding $\nu_i$ and $\alpha_i$ updates can be similarly modified as
\begin{align*}
\nu_i & \leftarrow \frac{\xi_i-\alpha_i}{|\xi_i-\alpha_i|}\frac{\Big[\sqrt{c_i}-|\xi_i-\alpha_i|\Big]_-}{\|\a_i\|^2},\\
\alpha_i & \leftarrow \frac{\xi_i-\alpha_i}{|\xi_i-\alpha_i|}\Big[\sqrt{c_i}-|\xi_i-\alpha_i|\Big]_-,
\end{align*}
and the rest of the updates in (\ref{algo:efficient}) stays the same. Conversely, if the constraint is a $\geq$ inequality, we only keep the nonnegative part of $\sqrt{c_i}-|\xi_i-\alpha_i|$ in the updates of $\nu_i$ and $\alpha_i$. If it is a bound constraint
\[
c_i - \epsilon \leq |\a_i^H\x|^2 \leq c_i + \epsilon,
\]
according to our previous discussion on ``rounding'' for this kind of constraint, we can define $\tau_i$ as
\[
\tau_i \leftarrow \begin{cases}
\sqrt{c_i-\epsilon}-|\xi_i-\alpha_i|, &
			\text{if } |\xi_i-\alpha_i|^2 < c_i - \epsilon, \\
\sqrt{c_i+\epsilon}-|\xi_i-\alpha_i|, &
			\text{if } |\xi_i-\alpha_i|^2 > c_i + \epsilon, \\
0,	& \text{otherwise.}
\end{cases}
\]
Then the corresponding updates of $\nu_i$ and $\alpha_i$ are
\begin{align*}
\nu_i & \leftarrow \frac{\xi_i-\alpha_i}{|\xi_i-\alpha_i|}\frac{\tau_i}{\|\a_i\|^2},\\
\alpha_i & \leftarrow \frac{\xi_i-\alpha_i}{|\xi_i-\alpha_i|}\tau_i.
\end{align*}

As we will see later, this type of memory efficient implementation can even be extended to cases when the constraints are not exactly homogeneous rank one quadratics. Furthermore, recall our previous discussion that a homogeneous rank one quadratic constraint is simply a linear constraint with an unknown phase (or sign in the real case), implying that if we have actual linear constraints in the QCQP problem, a similar idea can also be applied to avoid explicitly introducing a huge number of auxiliary variables, while still maintaining the simplicity of the updates.

\subsection{Initialization and parameter setting}\label{sec:phaseI}
At this point we need to remind the reader that, although the consensus ADMM algorithm we derived for non-convex QCQPs has an appealing form and cheap per-iteration complexity, it is after all a heuristic for what is in general an NP-hard problem. We may then anticipate that appropriate initialization and judicious parameter tuning will be more important than in standard applications of ADMM to convex problems. 
Nevertheless, we have devised practical rules that seem to work well in most cases, as discussed below.

The only parameter that needs to be tuned is $\rho$, which is only involved in the update of $\x$ if we have an explicit objective. Clearly a smaller $\rho$ steers the $\x$-update towards putting more emphasis on decreasing the cost function, whereas a bigger $\rho$ puts more weight on agreeing with the auxiliary variables $\{\z_i\}$, each guaranteed to satisfy one constraint. We found empirically that if we start with a feasible $\x$, then we can afford to have a relatively small value of $\rho$ for faster decrease of the cost, while preventing $\x$ from diverging towards infeasibility.

How can we find a feasible point for initialization? In some cases it is easy, for example when all the $\A_i$'s are positive semi-definite, and all the inequality constraints are homogeneous and one-sided, then a simple scaling suffices to make an arbitrary point feasible. In general, {\color{blue}finding a feasible point is also NP-hard}. In our context, we can attempt to find a feasible point by using the same consensus ADMM algorithm for the following feasibility problem, rewritten in the consensus form
\begin{equation*}
\begin{aligned}
\text{find}~~ & \x,\{\z_i\}_{i=1}^m, \\
\text{such that}~~	& \z_i^H\A_i\z_i - 2\real{\b_i^H\z_i} \leq c_i, \\
		& \z_i = \x,~\forall~i=1,...,m,
\end{aligned}
\end{equation*}
Applying consensus ADMM, we obtain the following updates
\begin{align*}
\x\; & \leftarrow \frac{1}{m}\sum_{i=1}^{m}(\z_i+\u_i), \\
\z_i\, & \leftarrow \argmin_{\z_i}~\|\z_i-\x+\u_i\|^2, \\
		& ~~~\st~~ \z_i^H\A_i\z_i - 2\real{\b_i^H\z_i} \leq c_i \\
\u_i & \leftarrow \u_i + \z_i - \x,
\end{align*}
which are completely independent of $\rho$ \footnote{Or one can interpret this as $\rho$ being $+\infty$.}. This type of iterates with random initialization 
usually converges much faster in finding a feasible point, if one exists. The result can then serve as initialization for subsequent ADMM updates with the cost function brought back into consideration. If the first phase fails to find a feasible point even after multiple trials, then we have, to some extent, numerical evidence that the problem may be infeasible, and for practical purposes there is often no point in proceeding further anyway.

\subsection{Parallel and distributed implementation}
Consensus ADMM is by its very nature highly parallelizable, since the update of each auxiliary variable $\z_i$ is independent of all others. This nice property is thankfully maintained even in the squeezed form (\ref{algo:efficient}), since all the operations involved are element-wise, except for two matrix vector multiplications $\A_s^H\x$ and $\A_s\bm{\nu}$, which can also be parallelized easily. This means that the proposed algorithm can easily achieve $p$-fold acceleration by using $p$ processors on a shared-memory system.

A more interesting case is when a large amount of data is stored in distributed storage, and different agents need to coordinate with a master node with small communication overheads. Suppose the data for the constraints $\{\A_i,\b_i\}$ are stored across $p$ agents, all connected to the master node which is in charge of the cost function. Since we assign each constraint an individual variable $\z_i$ and dual $\u_i$, suppose the $k$-th agent is in charge of $m_k$ constraints, a naive implementation would require the $j$-th agent to send $m_j$ of the $\z_i$'s and $\u_i$'s to the central node in each iteration. This is not necessary, as a matter of fact, since for the update of $\x$ only the sum of all the $\z_i$'s and the $\u_i$'s is required. Therefore, to minimize communication overheads, the $j$-th agent can simply define $\x_j = \sum_{i\in\varOmega_j}(\z_i+\u_i)$, where $\varOmega_j$ is the index set of the constraints handled by the $j$-th agent. At the master node, another summation over all the $\x_j$'s is carried out for the exact update of $\x$.


\section{Applications and numerical experiments}\label{sec:5}
So far we have introduced the general idea of applying consensus ADMM to QCQPs with efficient per-iteration updates, and explored memory-efficient and parallel/distributed implementation issues. In this section, we will look into some important QCQP applications, write down explicitly the algorithm, and compare its numerical performance with some state-of-the-art algorithms. All simulations were performed in MATLAB on a Linux desktop with 8 Intel i7 cores and 32GB of RAM.

\subsection{Feasible point pursuit}
One of the main drawbacks of the SDR approach for non-convex QCQPs is that when it is not obvious how to find a feasible point that satisfies the constraint set, there is a high chance that SDR, followed by taking the principal component and/or Gaussian randomization, will not satisfy all the constraints either. Recently, a new algorithm called FPP-SCA \cite{mehanna2015feasible} was proposed to address this issue by iteratively linearizing the non-convex part of the problem, while adding nonnegative slacks to each constraint and penalizing the sum of slacks in the cost function as well. Simulations in \cite{mehanna2015feasible} suggest that FPP-SCA works well with high probability, even when SDR fails.

Consensus ADMM can also be used to find feasible points, and it is possible to aim it towards finding a feasible point having smallest $\ell_2$ norm. Instead of giving each constraint a slack and trying to minimize the sum of the slacks, consensus ADMM gives each constraint a local variable and tries to drive these local variables to consensus. Explicitly, let us consider the following problem
\begin{equation}\label{prob:fpp}
\begin{aligned}
\minimize_{\x\in\C^n}~~ & \|\x\|^2 \\
\st~~ & \x^H\A_i\x \leq c_i,~\forall~i=1,...,m,
\end{aligned}
\end{equation}
where $\A_1,...,\A_m$ are in general Hermitian indefinite and full rank. Following our discussion in Sec.~\ref{sec:qcqp1-3}, the detailed consensus ADMM algorithm for (\ref{prob:fpp}) is given in Alg.~\ref{algo:fpp}, where we have applied the two stage approach described in Sec.~\ref{sec:phaseI}: we attempt to find a feasible point in the first phase, followed by stably decreasing its norm in the second phase. We found empirically that simply setting $\rho=1$ works very well for the second phase in this context.
\begin{algorithm}[t]
\LinesNumbered
initialize $\x$, $\z_i$ and $\u_i$\;
\For{$i=1,...,m$}{
Take the eigen-decomposition of $\A_i=\Q_i\L_i\Q_i^H$\;
}
\Repeat{$\x$ feasible}{
$\x \leftarrow \frac{1}{m}\sum_{i=1}^{m}\left(\z_i + \u_i\right)$\;
\For{$i=1,...m$}{
$\mu_i \leftarrow \arg_{\mu} \phi_i(\mu)=0$ using Alg.~1~or~2\;
$\mu_i \leftarrow \max\left\{0,\mu_i\right\}$\tcp*{for $\leq$ constraints}
$\z_i\;\! \leftarrow \Q_i\inv{\eye+\mu_i\L_i}\Q_i^H(\x-\u_i)$\;
$\u_i \leftarrow \u_i + \z_i - \x$\;
}}
$\rho=1$\;
\Repeat{The successive difference of $\x$ is smaller than $\varepsilon$}{
$\x \leftarrow \frac{1}{m + \rho^{-1}}\sum_{i=1}^{m}\left(\z_i + \u_i\right)$\;
\For{$i=1,...m$}{
$\mu_i \leftarrow \arg_{\mu} \phi_i(\mu)=0$ using Alg.~1~or~2\;
$\mu_i \leftarrow \max\left\{0,\mu_i\right\}$\tcp*{for $\leq$ constraints}
$\z_i\;\! \leftarrow \Q_i\inv{\eye+\mu_i\L_i}\Q_i^H(\x-\u_i)$\;
$\u_i \leftarrow \u_i + \z_i - \x$\;
}}
\caption{consensus-ADMM for (\ref{prob:fpp})}
\label{algo:fpp}
\end{algorithm}

Now let us compare consensus-ADMM with FPP-SCA on some synthetically generated problems. After fixing the problem dimension $n$ and $m$, we first generate $\x_\text{feas}\sim\CN(0,\eye)$. A Hermitian indefinite matrix $\A_i$ is generated by first randomly drawing a $n\times n$ matrix from $\CN(0,1)$, and then taking the average of its Hermitian and itself. The corresponding $c_i$ is set to be $\x_\text{feas}^H\A_i\x_\text{feas}^{} - |v_i|$ where $v_i$ is randomly generated from $\N(0,1)$. The constructed constraint set is therefore guaranteed to be non-empty, because we know $\x_\text{feas}$ is feasible, similar to the problem setting considered in \cite{mehanna2015feasible}. For $n=20$ and $m\in\{32,40,48\}$, the averaged results over 100 Monte-Carlo trials are presented in Table~\ref{tab:fpp}, and in each trial both ADMM and FPP-SCA are initialized with the same point, which is randomly generated from $\CN(0,\eye)$. As we can see, our proposed algorithm is able to produce similar performance with a much shorter execution time. It is possible to develop specialized solvers for FPP-SCA to accelerate it, but it is a non-trivial task which may require a lot of thinking, whereas our algorithm, readily available in Alg.~\ref{algo:fpp}, only requires elementary operations and simple iterations, thus it is also easy to code in a lower-level language.
\begin{table}[t]
\centering
\caption{Averaged performance, over 100 Monte-Carlo trials, in feasible point pursuit. In each column, consensus-ADMM is on the left, and FPP-SCA is on the right. The average loss is defined as $10\log_{10}(\|\x\|^2/\tr{\X})$, where $\X$ is the solution of the SDR.}
\label{tab:fpp}
\begin{tabular}{c|cc|cc|cc}
\hline
		& \multicolumn{2}{|c|}{feasible point} & \multicolumn{2}{|c|}{avg. loss (dB)}	
			& \multicolumn{2}{|c}{avg. time (sec.)}\\
\hline
$m=32$	& 100\%	& 100\%	& 0.376	& 0.375	& 4.5 & 31.9\\
$m=40$	& 100\%	& 100\%	& 0.503	& 0.526	& 5.4 & 37.7\\
$m=48$	& 100\%	& 100\%	& 0.600	& 0.597	& 7.9 & 44.5\\
\hline
\end{tabular}
\end{table}

\begin{figure}[!t]
\centering
\includegraphics[width=.8\textwidth]{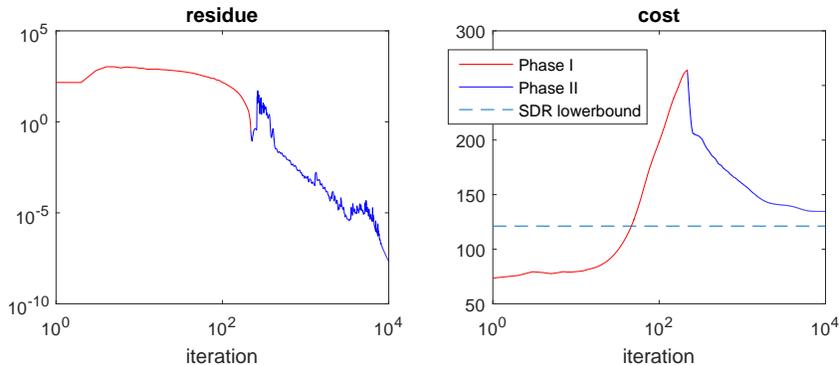}
\caption{Residual of the equality constraints $\sum_{i=1}^{m}\|\z_i-\x\|^2$ (left) and cost function $\|\x\|^2$ (right) vs. iteration number of Alg.~\ref{algo:fpp} for one random problem instance.}
\label{fig:fpp}
\end{figure}

To illustrate the scalability of our algorithm, we tested it on a larger problem with $n=100$ and $m=200$, and our algorithm took about 8 minutes to find a feasible point with smallest norm, which took about $10^4$ iterations. As shown in Fig.~\ref{fig:fpp}, the final result is not very far away from the generally unattainable lower bound provided by the SDR, with loss only about 0.45dB. If all we need is a feasible point, then it only requires about 200 iterations, showing great efficiency in finishing the most important task. In comparison, FPP-SCA requires more than 25 minutes to achieve a similar result.

\subsection{Multicast Beamforming}
Transmit beamforming is a wireless communication technique for transmitting signals to one or more users in a spatially selective way. A transmit beamforming system comprises a base station equipped with $n$ antennas, transmitting signals to a set of $m$ users within a certain service area, each having a single antenna. Assuming the transmitter is able to acquire channel state information (CSI) for all the users, multicast beamforming corresponds to the case when the base station \emph{selectively broadcasts} a common information-bearing signal to many users, ensuring a minimum received signal-to-noise ratio (SNR) at each user terminal, with the goal that the total transmission power is much smaller than the traditional method of radiating power isotropically around its service area. Clearly, this also controls interference to other nearby systems, e.g., in neighboring cells.

There are various formulations of multicast beamforming, ranging from single-group to multiple-groups, perfect channel state information (CSI) at the base station to channel second order statistics only, to name just a few; cf.~\cite{gershman2010convex} and the references therein. Almost all formulations are within the range of non-convex QCQP (and NP-hard~\cite{sidiropoulos2006transmit}), therefore it makes sense to test the performance of our proposed algorithm in this application. For brevity, we only consider the case when perfect CSI is available at the base station transmitter, corresponding to a fixed wireless scenario.

\subsubsection{Single-group multicast beamforming}
The most basic multicast beamforming formulation takes the following form~\cite{sidiropoulos2006transmit}
\begin{equation}\label{prob:mb1}
\begin{aligned}
\minimize_{\w\in\C^n}~~ & \|\w\|^2, \\
\st~~ & |\h_i^H\w|^2 \geq 1, \forall~i=1,...,m,
\end{aligned}
\end{equation}
where each $\h_i$ corresponds to the channel coefficients scaled according to the additive noise power and receive SNR requirement. Given $\{\h_i\}$, we wish to guarantee a certain SNR to all the receivers, while minimizing the transmit power $\|\w\|^2$.

Problem (\ref{prob:mb1}) is exactly in the form of (\ref{prob:rank_1}), except that the constraints are with inequalities, thus we can direly use the memory-efficient updates (\ref{algo:efficient}). For initialization, we found empirically that it is better to initialize with a feasible point, so that we can afford to use a smaller $\rho$ to provide faster convergence, reducing the risk of having unstable sequences. Fortunately it is easy to find a feasible point for (\ref{prob:mb1}) -- for any random initialization point, one only needs to scale it up until all the constraints are satisfied. The complete algorithm for (\ref{prob:mb1}), including the initialization strategy and our choice of $\rho$ is given in Alg.~\ref{algo:mb1}.

\begin{algorithm}[t]
\LinesNumbered
initialize $\w\sim\N(0,\eye)$\;
$\w \leftarrow \w / \min(|\bm{H}_s^H\w|)$\;
$\z_s \leftarrow m\w$\;
$\u_s \leftarrow 0$\;
$\rho=2\sqrt{m}$\;
\Repeat{The successive difference of $\w$ is smaller than $\varepsilon$}{
$\w \leftarrow \frac{1}{m+\rho^{-1}}\left(\z_s + \u_s\right)$\;
$\bm{\xi} \leftarrow \bm{H}_s^H\w$\;
$\nu_i \leftarrow \frac{\xi_i-\alpha_i}{|\xi_i-\alpha_i|}\frac{[1-|\xi_i-\alpha_i|]_+}{\|\h_i\|^2}$\;
$\z_s \leftarrow m\w - \u_s + \bm{H}_s\bm{\nu}$\;
$\u_s \leftarrow \u_s + \z_s - m\w$\;
$\alpha_i \leftarrow \frac{\xi_i-\alpha_i}{|\xi_i-\alpha_i|}\Big[1-|\xi_i-\alpha_i|\Big]_+$\;
}
\caption{consensus-ADMM for (\ref{prob:mb1})}
\label{algo:mb1}
\end{algorithm}

\begin{figure}[t!]
\centering
\includegraphics[width=.8\textwidth]{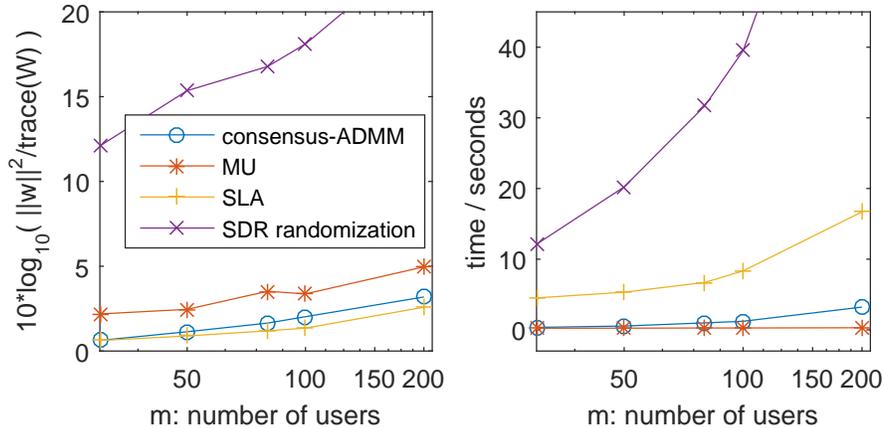}
\caption{Averaged performance of various methods for (\ref{prob:mb1}) over 100 Monte-Carlo trials, with $n=100$ and $m$ from $30$ to $200$. Performance gap relative to the generally unattainable SDR lower bound on the left, and computation time on the right.}
\label{fig:mb1}
\end{figure}

We test the numerical performance of Alg.~\ref{algo:mb1} on problem (\ref{prob:mb1}), and we compare it with the standard SDR followed by Gaussian randomization~\cite{sidiropoulos2006transmit}, successive linear approximation (SLA)~\cite{tran2014conic}, and the recently proposed multiplicative update (MU)~\cite{gopalakrishnan2015high}, which uses an approximate formulation and leads to highly efficient approximate solutions to (\ref{prob:mb1}). CVX~\cite{cvx} is used to solve the SDR and SLA in our experiment. We fix $n=100$, and vary $m \in \{30,50,80,100,200\}$, with each channel vector $\h_i$ randomly generated from $\CN(0,\eye)$. The averaged gap to the SDR lower bound, and the averaged computation time for all algorithms, with averages taken over 100 Monte-Carlo trials, are shown in Fig.~\ref{fig:mb1}, with each algorithm initialized at the same randomly generated points. As we can see, ADMM is able to give slightly worse performance than SLA in terms of transmission power, while requiring far smaller execution time.

To test the scalability of our algorithm, we also applied it to a massive MIMO multicast scenario with $n=500$ antennas and $m=100$ users. In this case, CVX is not able to solve the SDR within reasonable amount of time, so we can only compare the transmission power $\|\w\|^2$ without knowing how far it is from the SDR lower bound. It is reported in~\cite{gopalakrishnan2015high} that MU followed by one step of SLA gives the best result in both minimizing the transmission power and keeping computation time low, so we compare our algorithm with this two-step strategy here, with one or up to ten SLA iterations (unless $\|\w\|^2$ converges, i.e., per-iteration improvement is less than $10^{-5}$). For fair comparison, ADMM is initialized with the same point generated by MU. ADMM is able to compute a beamformer that is better than one step of SLA and do so in less time; in fact the transmission power obtained via ADMM is only slightly worse than ten steps of SLA. Notice that the update rule of ADMM is as simple as that of MU; both are simple enough to implement in real communication hardware, whereas SLA requires a full-blown convex optimization solver -- which seems unrealistic for base station deployment as of this writing.

\begin{table}[t!]
\centering
\caption{Averaged performance of various methods for (\ref{prob:mb1}) over 100 Monte-Carlo trials, with $n=500$ and $m=100$. }
\label{tab:mb1}
\begin{tabular}{c|c|c}
\hline
	& ~~~~~$\|\w\|^2$~~~~~ & computation time \\
\hline
consensus-ADMM & 0.1131 & 0.5235 sec. \\
1 step of SLA & 0.1213 & 0.6204 sec. \\
$\leq 10$ steps of SLA & 0.1125 & 6.2115 sec. \\
\hline
\end{tabular}
\end{table}

\subsubsection{Secondary user multicast beamforming}
We now consider adding primary user interference constraints to the basic multicast beamforming formulation in (\ref{prob:mb1}). This secondary multicast underlay scenario has been considered in \cite{phan2009spectrum}, and the problem of interest can be formulated as
\begin{equation}\label{prob:mb2}
\begin{aligned}
\minimize_{\w\in\C^n}~~ & \|\w\|^2, \\
\st~~ & |\h_i^H\w|^2 \geq \tau, \forall~i=1,...,m,\\
	  & |\g_k^H\w|^2 \leq \eta, \forall~k=1,...,l,\\
\end{aligned}
\end{equation}
where in this case we have $l$ additional primary users who should be protected for excess interference cause by the secondary multicast transmission, and $\g_k$ denotes the channel vector from the multicast transmitter to the $k$-th (single-antenna) primary user receiver.

Again, the efficient updates in (\ref{algo:rank_1}) for rank one quadratic constraints can be used. As for initialization, there is no obvious way to find a feasible point in this case, so the two-stage procedure we discussed before is used, which ignores the cost function first to find a feasible point, and then uses this feasible point to initialize the complete updates with a relatively small $\rho$ to accelerate convergence. The complete algorithm is given in Alg.~\ref{algo:mb2}.

\begin{algorithm}[t]
\LinesNumbered
initialize $\w\sim\N(0,\eye)$\;
$\z_s \leftarrow (m+l)\w$\;
$\u_s \leftarrow 0$\;
\Repeat{$\w$ feasible}{
$\w \leftarrow \frac{1}{m+l}\left(\z_s + \u_s\right)$\;
$\begin{array}{l}
\bm{\xi} \leftarrow \bm{H}_s^H\w \\
\bar{\bm{\xi}}\leftarrow \bm{G}_s^H\w
\end{array}$\;
$\begin{array}{l}
\nu_i \leftarrow \frac{\xi_i-\alpha_i}{|\xi_i-\alpha_i|}
\frac{[\sqrt{\tau}-|\xi_i-\alpha_i|]_+}{\|\h_i\|^2} \\
\bar{\nu}_i \leftarrow \frac{\bar{\xi}_i-\bar{\alpha}_i}{|\bar{\xi}_i-\bar{\alpha}_i|}
\frac{[\sqrt{\eta}-|\bar{\xi}_i-\bar{\alpha}_i|]_-}{\|\g_i\|^2}
\end{array} $\;
$\z_s \leftarrow (m+l)\w - \u_s + \bm{H}_s\bm{\nu} + \bm{G}_s\bar{\bm{\nu}}$\;
$\u_s \leftarrow \u_s + \z_s - (m+l)\w$\;
$\begin{array}{l}
\alpha_i \leftarrow \frac{\xi_i-\alpha_i}{|\xi_i-\alpha_i|}\big[\sqrt{\tau}-|\xi_i-\alpha_i|\big]_+ \\
\bar{\alpha}_i \leftarrow \frac{\bar{\xi}_i-\bar{\alpha}_i}{|\bar{\xi}_i-\bar{\alpha}_i|}
\big[\sqrt{\eta}-|\bar{\xi}_i-\bar{\alpha}_i|\big]_-
\end{array}$\;
}
$\rho=2\sqrt{m+l}$\;
\Repeat{The successive difference of $\w$ is smaller than $\varepsilon$}{
$\w \leftarrow \frac{1}{m+l+\rho^{-1}}\left(\z_s + \u_s\right)$\;
$\begin{array}{l}
\bm{\xi} \leftarrow \bm{H}_s^H\w \\
\bar{\bm{\xi}}\leftarrow \bm{G}_s^H\w
\end{array}$\;
$\begin{array}{l}
\nu_i \leftarrow \frac{\xi_i-\alpha_i}{|\xi_i-\alpha_i|}
\frac{[\sqrt{\tau}-|\xi_i-\alpha_i|]_+}{\|\h_i\|^2} \\
\bar{\nu}_i \leftarrow \frac{\bar{\xi}_i-\bar{\alpha}_i}{|\bar{\xi}_i-\bar{\alpha}_i|}
\frac{[\sqrt{\eta}-|\bar{\xi}_i-\bar{\alpha}_i|]_-}{\|\g_i\|^2}
\end{array} $\;
$\z_s \leftarrow (m+l)\w - \u_s + \bm{H}_s\bm{\nu} + \bm{G}_s\bar{\bm{\nu}}$\;
$\u_s \leftarrow \u_s + \z_s - (m+l)\w$\;
$\begin{array}{l}
\alpha_i \leftarrow \frac{\xi_i-\alpha_i}{|\xi_i-\alpha_i|}\big[\sqrt{\tau}-|\xi_i-\alpha_i|\big]_+ \\
\bar{\alpha}_i \leftarrow \frac{\bar{\xi}_i-\bar{\alpha}_i}{|\bar{\xi}_i-\bar{\alpha}_i|}
\big[\sqrt{\eta}-|\bar{\xi}_i-\bar{\alpha}_i|\big]_-
\end{array}$\;
}
\caption{consensus-ADMM for (\ref{prob:mb2})}
\label{algo:mb2}
\end{algorithm}

\begin{figure}[t!]
\centering
\includegraphics[width=.8\textwidth]{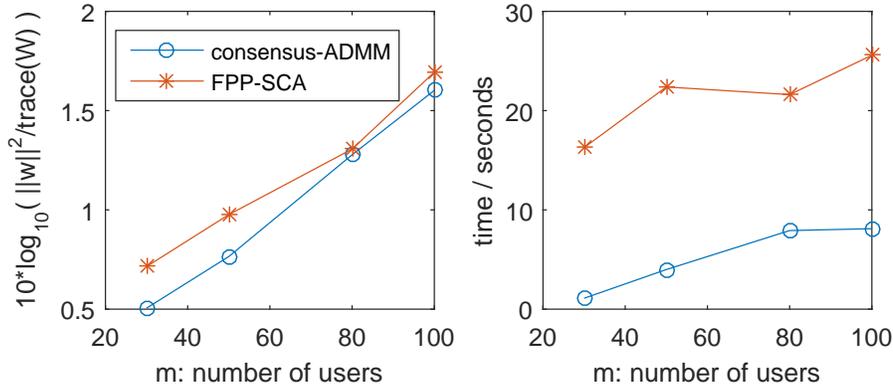}
\caption{Averaged performance of ADMM and FPP-SCA for (\ref{prob:mb2}) over 100 Monte-Carlo trials, with performance degradation relative to the generally unattainable SDR lower bound on the left, and computation time on the right.}
\label{fig:mb2}
\end{figure}

Similar to the previous simulation settings, we fix $n=100$ and $l=10$, and vary $m\in\{30,50,80,100\}$, with channel coefficients randomly generated from $\CN(0,\eye)$. For $\tau=10$ and $\eta=1$, the averaged performance over 100 Monte-Carlo runs comparing to the SDR lower bound and FPP-SCA as described in~\cite{mehanna2015feasible} is shown in Fig. \ref{fig:mb2}. SDR randomization almost always fails to find a feasible point in this case, thus not considered in this experiment. Both methods are initialized with a random point from $\CN(0,\eye)$, and manage to obtain a good feasible point in all problem instances that we tried. We allow ADMM to take multiple initializations if the first stage of Alg.~\ref{algo:mb2} does not find a feasible point after $10^3$ iterations, thus the computation time of ADMM is more likely to vary (as seen on the right panel of Fig.~\ref{fig:mb2} for $m=80$ and $100$), although much smaller than that of FPP-SCA, which requires a general-purpose convex optimization solver, in our case CVX~\cite{cvx}. Note that ADMM also yields better performance than FPP-SCA in terms of transmit power $\|\w\|^2$.

\subsection{Phase retrieval}
Phase retrieval is the problem of estimating a signal from the magnitude of complex linear measurements, without access to the corresponding phases. This problem arises in various applications like crystallography, microscopy, and optical imaging \cite{shechtman2015phase}. Specifically, let $\s$ be the desired signal, the measurements $\{y_i\}_{i=1}^m$ are collected via $y_i=|\a_i^H\s|^2$, possibly perturbed by noise. In the sequel we will see that for a number of noise models this problem can be written as non-convex QCQP, therefore we can test our algorithm together with other state-of-the-art phase retrieval methods. Notice that for some specific measurement systems the problem actually has hidden convexity, e.g., this is the case for 1-D over-sampled Fourier measurements~\cite{huang2016icassp}, but here we focus on a general measurement setup.

\subsubsection{Noiseless case}
Assuming all the measurements are exact, we can write the phase retrieval problem as the following feasibility problem
\begin{equation}\label{prob:pr1}
\begin{aligned}
\text{find~~} & \x \in \C^n, \\
\text{such that~~} & |\a_i^H\x|^2 = y_i,~ \forall~i=1,...,m.
\end{aligned}
\end{equation}
This is exactly in the form of (\ref{prob:rank_1}) except that there is no explicit cost function, so we can apply the memory efficient implementation of consensus-ADMM to obtain the following updates
\begin{equation}\label{alg:pr1}
\boxed{
\begin{aligned}
\x   & \leftarrow \frac{1}{m}\left(\z_s + \u_s\right), \\
\bm{\xi} & \leftarrow \A_s^H\x, \\
\nu_i  & \leftarrow \frac{\xi_i-\alpha_i}{|\xi_i-\alpha_i|}\frac{\sqrt{y_i}-|\xi_i-\alpha_i|}{\|\a_i\|^2}, \\
\z_s & \leftarrow m\x - \u_s + \A_s\bm{\nu}, \\
\u_s & \leftarrow \u_s + \z_s - m\x, \\
\alpha_i & \leftarrow \frac{\xi_i-\alpha_i}{|\xi_i-\alpha_i|}\left(\sqrt{y_i}-|\xi_i-\alpha_i|\right),
\end{aligned}
}
\end{equation}
where $\A_s = [~\a_1~\a_2~...~\a_m~]$ is obtained by stacking all the $\a_i$ vectors as its columns. Notice that since we do not have an explicit cost function, it does not matter what value of $\rho$ we choose -- they all work the same for this problem.


\subsubsection{Bounded noise}
In practice the measurements are seldom perfect, so we need to incorporate uncertainties in the measurements. A simple assumption is that measurements are quantized at relatively high resolution, in which case we can model the measurements as being corrupted by noise that is uniformly distributed between $[-\epsilon,\epsilon]$. We can modify the noiseless formulation (\ref{prob:pr1}) as follows (similar to the B-FPP formulation proposed in~\cite{qian2015phase})
\begin{equation}\label{prob:pr2}
\begin{aligned}
\text{find~~} & \x \in \C^n, \\
\text{such that~~} & y_i-\epsilon \leq |\a_i^H\x|^2 \leq y_i+\epsilon,~ \forall~i=1,...,m.
\end{aligned}
\end{equation}
Again we can apply the memory efficient implementation with the following updates
\begin{equation}\label{alg:pr2}
\boxed{
\begin{aligned}
\x   & \leftarrow \frac{1}{m}\left(\z_s + \u_s\right), \\
\bm{\xi} & \leftarrow \A_s^H\x, \\
\tau_i &\leftarrow {\small \begin{cases}
\sqrt{y_i-\epsilon}-|\xi_i-\alpha_i|, &
			\text{if } |\xi_i-\alpha_i|^2 < y_i - \epsilon, \\
\sqrt{y_i+\epsilon}-|\xi_i-\alpha_i|, &
			\text{if } |\xi_i-\alpha_i|^2 > y_i + \epsilon, \\
0,	& \text{otherwise,}
\end{cases} }\\
\nu_i & \leftarrow \frac{\xi_i-\alpha_i}{|\xi_i-\alpha_i|}\frac{\tau_i}{\|\a_i\|^2},\\
\z_s & \leftarrow m\x - \u_s + \A_s\bm{\nu}, \\
\u_s & \leftarrow \u_s + \z_s - m\x, \\
\alpha_i & \leftarrow \frac{\xi_i-\alpha_i}{|\xi_i-\alpha_i|}\tau_i.
\end{aligned}
}
\end{equation}

\subsubsection{Gaussian noise}
Another interesting scenario is that where measurements are corrupted by additive white Gaussian noise, in which case maximum likelihood estimation can be cast as the following non-convex QCQP (similar to the LS-FPP formulation proposed in~\cite{qian2015phase})
\begin{equation}\label{prob:pr3}
\begin{aligned}
\minimize_{\x\in\C^n,\w\in\R^m}~~ & \frac{1}{2}\|\w\|^2 \\
\st~~ & |\a_i^H\x|^2 = y_i + w_i,~ \forall~i=1,...,m.
\end{aligned}
\end{equation}
This kind of constraint is not covered in our previous discussions, so we study this case in a bit more detail here. Let us first rewrite (\ref{prob:pr3}) into a consensus optimization form by introducing $m$ auxiliary variables $\z_1,...,\z_m$ replicating $\x$ for each constraint
\begin{align*}
\minimize_{\x,\{\z_i\},\w}~~ & \half\|\w\|^2 \\
\st~~ & |\a_i^H\z_i|^2 = y_i + w_i, \\
		& \z_i = \x, ~ \forall~i=1,...,m.
\end{align*}
The plain vanilla version of ADMM, treating $\x$ as the first block and $\{\z_i\}$ and $\w$ as the second block, takes the following form
\begin{align*}
\x &\leftarrow \frac{1}{m}\sum_{i=1}^{m}\left(\z_i+\u_i\right), \\
(\z_i,w_i) &\leftarrow \argmin_{|\a_i^H\x|^2 = y_i + w_i} \half|w_i|^2 + \rho\|\z_i-\x+\u_i\|^2, \\
\u_i &\leftarrow \u_i+\z_i-\x.
\end{align*}
The main difficulty boils down to an efficient method for the second update, which can be written explicitly as
\begin{equation}\label{prob:pr3z}
\begin{aligned}
\minimize_{\z_i,w_i}~~ & \half|w_i|^2 + \rho\|\z_i-\x+\u_i\|^2 \\
\st~~ & |\a_i^H\z_i|^2 = y_i + w_i.
\end{aligned}
\end{equation}

The same idea of using the Lagrangian can be applied to solve (\ref{prob:pr3z}). Constructing the Lagrangian with a single multiplier $\mu_i$ and setting its derivative with respect to $w_i$ and $\z_i$ equal to 0, we have
\[
w_i = \mu_i,
\]
which, interestingly, means that the optimal multiplier is actually equal to the estimated noise term, and
{\color{blue}
\begin{equation}\label{eq:pr3z}
\begin{aligned}
\z_i &= \inv{\rho\eye+\mu_i\a_i\a_i^H}\rho(\x-\u_i) \\
	&= \x-\u_i-\frac{\mu_i\a_i^H(\x-\u_i)}{\rho+\mu_i\|\a_i\|^2}\a_i.
\end{aligned}
\end{equation}}
Plugging them back into the equality constraint, we end up with an equation with respect to $\mu_i$
{\color{blue}\begin{equation}\label{eq:pr3mu}
\frac{\rho^2|\a_i^H(\x-\u_i)|^2}{(\rho+\|\a_i\|^2\mu_i)^2} = y_i + \mu_i.
\end{equation}}
Equation (\ref{eq:pr3mu}) is equivalent to a cubic equation, for which we know the formula for the three roots. Moreover, since we know the three roots of a real cubic equation are either all real or one real and two complex conjugates, and that the correct $\mu_i$ we are looking for is real, we can deliberately select the value of $\rho$ so that the latter case happens, resolving the ambiguity in solving (\ref{eq:pr3mu}). Detailed derivation of the formula for solving (\ref{eq:pr3mu}) is given in Appendix~\ref{Appendix:A}, where it is also shown that by setting $\rho = 1.1\max_i y_i\|\a_i\|^2$, each equation (\ref{eq:pr3mu}) is guaranteed to have a unique real root, thus being the correct multiplier we are looking for.

Memory efficient implementation is again applicable here, which eventually leads to the following updates:
\begin{equation}\label{alg:pr3}
\boxed{
\begin{aligned}
\x   & \leftarrow \frac{1}{m}\left(\z_s + \u_s\right), \\
\bm{\xi} & \leftarrow \A_s^H\x, \\
\mu_i &\leftarrow \text{steps described in Appendix~\ref{Appendix:A}}, \\
\nu_i & \leftarrow -\frac{\mu_i}{\rho+\mu_i\|\a_i\|^2}\left(\xi_i-\alpha_i\right), \\
\z_s & \leftarrow m\x-\u_s+\A_s\bm{\nu},\\
\u_s & \leftarrow \u_s+\z_s-m\x,\\
\alpha_i &\leftarrow -\frac{\mu_i\|\a_i\|^2}{\rho+\mu_i\|\a_i\|^2}\left(\xi_i-\alpha_i\right).
\end{aligned}
}
\end{equation}
The only unclear part in (\ref{alg:pr3}) is the update of $\mu_i$. However, since $\mu_i$ is a solution of (\ref{eq:pr3mu}), which only depends on $\a_i^H(\x+\u_i)=\xi_i+\alpha_i$, it is indeed possible to implement (\ref{alg:pr3}) without explicitly calculating the individual $\z_i$'s and $\u_i$'s.

\subsubsection{Adding priors}
In a lot of cases there is useful prior information available about the signal that can help enhance the estimation performance. For example, one may know a priori that the signal to be estimated is real, non-negative, and/or sparse. All of these type of prior information can easily be incorporated into the $\x$-update, which usually boils down to very simple projections, like zeroing out the imaginary part and/or zeroing out the negative values.

For sparsity, a popular method is to add an $\ell_1$ penalty to the cost, which in our case leads to a {\em soft-thresholding} to $\frac{1}{m}\left(\z_s + \u_s\right)$ for the update of $\x$. However, recall that we are dealing with NP-hard problems and there is no guarantee that our algorithm is always going to reach a global solution. Therefore, for practical purposes it is sometimes better to just use the straight-forward cardinality constraint, which is easy to tune for the desired sparsity level, and has an equally simple {\em hard-thresholding} update.

{\color{blue}
\subsubsection{Simulations}
Now we show some numerical results on the performance of the proposed algorithms. There exist many algorithms designed for phase retrieval under specific scenarios, for example, the classical Gerchberg-Saxton algorithm~\cite{gerchberg1972practical} and Fienup's algorithm~\cite{fienup1978reconstruction} were designed for Fourier-based measurements, and have been successfully applied in the phase retrieval community for decades. More recently, inspired by the success of compressive sensing, random Gaussian measurements have gained more attention, and the state-of-the-art algorithms include alternating minimization~\cite{netrapalli2015phase} and Wirtinger flow~\cite{candes2015phase}, both of which contain a special initialization step to help convergence. SDR-based methods have also been developed, including PhaseLift~\cite{candes2013phase} and PhaseCut~\cite{waldspurger2015phase}, however, they again suffer from effectively squaring the problem dimension, and how to recover a good approximate solution to the original problem when they return a higher rank matrix is an open question.

We consider the columns of the sensing matrix $\{\a_i\}$ to be generated from an i.i.d. complex Gaussian distribution $\CN(0,\eye)$. In a noiseless scenario, we apply the updates (\ref{alg:pr1}) to the problem of phase retrieval with random Gaussian measurements, and compare it with two state-of-the-art algorithms designed for this setting: Wirtinger flow~\cite{candes2015phase} and alternating minimization~\cite{netrapalli2015phase}. Notice that algorithm (\ref{alg:pr1}) only requires two matrix-vector multiplications, thus the per-iteration complexity is the same as that of Wirtinger flow and alternating minimization.

We randomly generate a desired signal $\s\in\C^n$ where $n=128$, then take $m$ phase-less measurements $|\a_i^H\s|^2$ with $m$ ranging from $2n$ to $5n$. Using the same initialization proposed in~\cite{netrapalli2015phase} and~\cite{candes2015phase}, we let consensus-ADMM (Algorithm~\eqref{alg:pr1}), Wirtinger flow, and alternating minimization run for at moat $10^5$ iterations, and the probability of resolution over 100 Monte-Carlo trials is given in Fig.~\ref{fig:pr_gm}, where we declare that the signal has been resolved if (after adjusting for the global phase ambiguity of the estimate $\x$),
\[
\min_{\theta}\|e^{j\theta}\x-\s\|^2 < 10^{-5}.
\]
It is very satisfying to see that consensus-ADMM has higher empirical probability of resolution than both Wirtinger flow and alternating minimization. Due to the similar per-iteration complexities of the three algorithms, running $10^5$ iterations take approximately the same time for all of them.
\begin{figure}[t!]
\centering
\includegraphics[width=.8\textwidth]{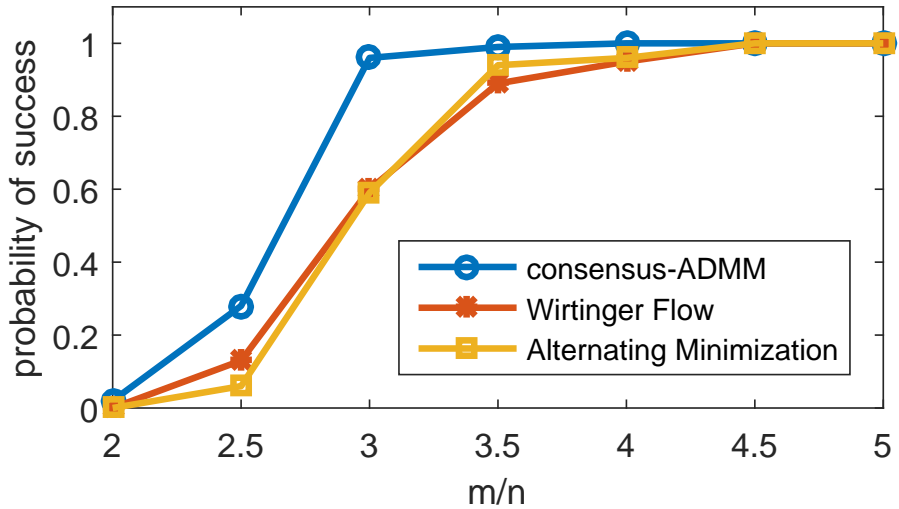}
\caption{Empirical probability of resolution based on 100 Monte-Carlo trials for various number of measurements.}
\label{fig:pr_gm}
\end{figure}

Finally, we briefly show the performance of Algorithm~\eqref{alg:pr2} and~\eqref{alg:pr3}, under their corresponding noise models, and we only show the results for $n=128$ and $m=5n$, with $\s$ and $\A_s$ generated as before. Consider the quantized measurements $\y = \lfloor|\A_s^H\s|^2\rceil$, where $\lfloor\cdot\rceil$ rounds the argument to the nearest integer, we can use formulation \eqref{prob:pr2} with $\epsilon=0.5$, and apply Algorithm~\eqref{alg:pr2}. The number of constraint violations and the mean squared error (MSE) are given in Table~\ref{tab:pr2}, each averaged over 100 Monte Carlo trials, where MSE is defined as
\[
\text{MSE}=10\log_{10}\left(\min_{\theta}\|e^{j\theta}\x-\s\|^2\right).
\]
As we can see, Algorithm~\eqref{alg:pr2} is able to give a solution that is consistent with all the measurements in all cases, whereas the other two algorithms cannot, even though their MSE performance is still pretty good. For additive white Gaussian noise $\y=|\A_s^H\s|^2 +\w$, we fix the SNR to be 20dB, and the averaged performance over 100 Monte Carlo trials is shown in Table~\ref{tab:pr3}. Algorithm~\eqref{alg:pr3} performs almost as well as Wirtinger flow, and both perform better than alternating minimization in this case. This is as expected, since alternating minimization aims to solve a different formulation, which is not the maximum likelihood one for this model. To sum up, consensus-ADMM is able to achieve similar (if not better) performance to the state-of-the-art methods for phase retrieval with random Gaussian measurements.

\begin{table}[t!]
\centering
\caption{Performance of quantized phase retrieval}
\label{tab:pr2}
\begin{tabular}{c|c|c}
\hline
 & \# of violations & MSE \\
\hline
consensus-ADMM & 0 & -37dB \\
Wirtinger flow & 13.3 & -31dB \\
Alternating Min. & 143.4 & -34dB \\
\hline
\end{tabular}
\end{table}

\begin{table}[t!]
\centering
\caption{Performance of Gaussian noise phase retrieval}
\label{tab:pr3}
\begin{tabular}{c|c|c}
\hline
 & & \\[-7pt]
 & $\|\y-|\A_s\x|^2\|^2$ & MSE \\
\hline
consensus-ADMM & 6.9e3 & -2.4dB \\
Wirtinger flow & 6.9e3 & -2.4dB\\
Alternating Min. & 1.1e4 & 0.4dB\\
\hline
\end{tabular}
\end{table}
}

\section{Conclusion}\label{sec:6}
In this paper, we have proposed a new algorithm for general non-convex QCQPs, which is very different from any existing methods, general or specialized, for such problems. The main ideas behind this proposed algorithm are:
\begin{itemize}
\item Any QCQP-1 can be optimally solved, irrespective of (non-)convexity;
\item Consensus ADMM can be used to solve general QCQPs, in such a way that each update requires to solve a number of QCQP-1's.
\end{itemize}
For the type of QCQP-1's encountered here, we showed that they can be solved very efficiently, and made important steps towards enhancing scalability of the overall algorithm, including
\begin{itemize}
\item Memory-efficient implementations for certain types of ``simple'' constraints, e.g., rank one quadratics;
\item Parallel/distributed implementations with small communication overhead.
\end{itemize}

The proposed algorithm and various custom implementations were fully fleshed out and applied to various important non-convex QCQP problems, from pure feasibility pursuit to two real-world engineering tasks: multicast beamforming and phase retrieval. the proposed algorithms consistently exhibited favorable performance compared to the prior state of art, including classical and more modern methods for phase retrieval, which has drawn renewed interest in recent years. We believe that the general applicability of the proposed algorithm has the potential to bring significant performance improvements to many other applications of non-convex QCQP as well.

\appendix
\section{Proof of Theorem~\ref{thm:convergence}}\label{Appendix:proof}
A KKT point $\x_\star$ of (\ref{prob:qcqp}), together with the corresponding dual variable $\bm{\mu}_\star$, satisfies that
\begin{subequations}\label{eq:kkt}
\begin{align}
\A_0\x_\star - \b_0 + \sum_{i=1}^{m}\mu_{i\star}\left( \A_i\x_\star-\b_i \right) = 0, \\
\bm{\mu}_\star \geq 0, \\
\x_\star^H\A_i\x_\star - 2\real{\b_i^H\x_\star} \leq c_i, \\
\mu_{i\star}\left( \x_\star^H\A_i\x_\star - 2\real{\b_i^H\x_\star} \right) = 0, \\
\forall~i=1,...,m. \nonumber
\end{align}
\end{subequations}

{\color{blue}
Let us use a superscript $t$ to denote the point obtained at iteration $t$ by the update rule (\ref{alg:admm}), then at iteration $t+1$, since we assume that each $\{\z_i^t\}$ is well defined, we have that
\begin{subequations}\label{eq:update}
\begin{align}
\A_0\x^{t+1} + m\rho\x^{t+1} = \b_0 + \rho\sum_{i=1}^{m}(\z_i^t + \u_i^t), \\
\exists \mu_i^{t+1}\geq 0~\text{s.t.~}
\u_i^{t+1} + \mu_i^{t+1}\left( \A_i\z_i^{t+1}-\b_i \right) = 0,\\
\forall~i=1,...,m \nonumber
\end{align}
\end{subequations}
where in the second equation we used the fact that
\[
\u^{t+1}=\u^t+\z_i^{t+1}-\x^{t+1}.
\]

Now by assuming $\z_i^{t+1}-\x^{t+1} \rightarrow 0$, we trivially have that
\[
\u_i^{t+1}-\u_i^t \rightarrow 0.
\]
Therefore, for $t$ sufficiently large, we have
\[
\A_0\x^{t+1} + m\rho\x^{t+1} = \b_0 + \rho\sum_{i=1}^{m}\left(\x^t - \mu_i^{t+1}\left( \A_i\x^{t+1}-\b_i \right)\right).
\]
By further assuming
\[
\x^{t+1}-\x^t \rightarrow 0,
\]
this becomes exactly (\ref{eq:kkt}a) by setting $\bm{\mu}_\star=\rho\bm{\mu}^{t+1}$. The rest of the KKT conditions are guaranteed by the feasibility of $\z_i^{t+1}$ for the $i$-th constraint, and our assumption that $\z^{t+1}-\x^{t+1}\rightarrow0$.
}

\section{Solving (\ref{eq:pr3mu})}\label{Appendix:A}
We derive the solution for equation (\ref{eq:pr3mu}), which can be equivalently written as the following cubic equation after we drop the subscripts,
\begin{equation}\label{eq:cubic}
\begin{aligned}
\|\a\|^4\mu^3 + (2\rho\|\a\|^2+y\|\a\|^4)\mu^2 + (2y\rho\|\a\|^2+\rho^2)\mu \\
 + y\rho^2 - \rho^2|\a^H(\x-\u)|^2 = 0.
\end{aligned}
\end{equation}

For a general cubic equation
\[
\gamma_3\mu^3 + \gamma_2\mu^2 + \gamma_1\mu + \gamma_0 = 0,
\]
the three roots can be found with the following formulas (assuming all of these quantities are non-zero, which can be ensured by our specific choice of $\rho$ presented in the sequel)
\begin{align*}
\Delta_0 &= \gamma_2^2 - 3\gamma_3\gamma_1, \\
\Delta_1 &= 2\gamma_2^3 - 9\gamma_3\gamma_2\gamma_1 + 27\gamma_3^2\gamma_0, \\
C~ &= \sqrt[3]{\frac{\Delta_1 + \sqrt{\Delta_1^2-4\Delta_0^3}}{2}}, \\
\hat{\mu}_k &= -\frac{1}{3\gamma_3}\left(\gamma_2+\iota_kC+\frac{\Delta_0}{\iota_kC}\right),
\end{align*}
where $\iota_1 = 1, \iota_2 = \frac{-1+j\sqrt{3}}{2}, \iota_3 = \frac{-1-j\sqrt{3}}{2}$ are the three cubic roots of $1$.
Furthermore, if all the coefficients are real, then there is at least one real root, and the other two are either complex conjugates or both real, depending on whether $\Delta_1^2-4\Delta_0^3$ is negative or positive.

Now, let us plug in the coefficients of (\ref{eq:cubic}) into the formula, to get
\begin{align*}
\Delta_0 &= (\rho\|\a\|^2 - y\|\a\|^4)^2, \\
\Delta_1 &= -2(\rho\|\a\|^2 - y\|\a\|^4)^3 - 27\rho^2\|\a\|^8|\a^H(\x+\u)|^2, \\
\Delta_1^2-&4\Delta_0^3 = \Big(27\rho^2\|\a\|^8|\a^H(\x+\u)|^2\Big)^2 + \\
&\Big(27\rho^2\|\a\|^8|\a^H(\x+\u)|^2\Big) \Big(4(\rho\|\a\|^2 - y\|\a\|^4)^3\Big).
\end{align*}
Before we proceed, recall that the coefficients of (\ref{eq:pr3mu}) are all real, and the root that we are looking for, which is the optimal Lagrange multiplier, is also real. Therefore, to make our life easier, we can choose the value of $\rho$ to ensure that $\Delta_1^2-4\Delta_0^3 > 0$, for example by setting $\rho>y\|\a\|^2$. Then we can proceed to the rest of the steps and choose the solution to be $\hat{\mu}_1$, the only real root of (\ref{eq:pr3mu}). Putting the subscripts back and considering there are $m$ of them, this means we should set $\rho > \max_i y_i\|\a_i\|^2$. In practice smaller $\rho$ usually leads to faster convergence, we therefore recommend setting
\[
\rho = 1.1\max_i y_i\|\a_i\|^2.
\]

\section*{Acknowledgment}
The authors would like to thank Prof. Veit Elser from Cornell University for correcting a mistake made in an earlier version of this paper.

\bibliographystyle{unsrt}
\bibliography{refs}
\end{document}